\documentclass[12pt]{amsart}
\usepackage{amscd, amssymb}
\input xy
\xyoption{all}

\theoremstyle{plain}
\newtheorem{theorem}{Theorem}
\newtheorem{proposition}{Proposition}
\newtheorem{lemma}{Lemma}
\newtheorem{corollary}{Corollary}
\newtheorem*{lemma*}{Lemma}
\newtheorem*{theorem*}{Theorem}
\newtheorem*{proposition*}{Proposition}
\newtheorem*{corollary*}{Corollary}
\newtheorem*{conjecture*}{Conjecture}

\theoremstyle{definition}

\newtheorem*{definition*}{Definition}

\theoremstyle{remark}
\newtheorem{remark}{Remark}

\newtheorem*{convention}{Convention}
\newtheorem*{acknowledgements}{Acknowledgements}

\newcommand{\pp}{\mathbb{P}}
\newcommand{\zz}{\mathbb{Z}}
\newcommand{\qq}{\mathbb{Q}}
\newcommand{\cc}{\mathbb{C}}

\newcommand{\mbar}{{\overline{M}}}

\newcommand{\so}{{\mathcal{O}}}

\newcommand{\cl}{\mathcal{L}}
\newcommand{\tlb}{\underline{1}}
\newcommand{\h}{\mathcal{H}}

\newcommand{\p}{\partial}
\newcommand{\lp}{\left(}
\newcommand{\rp}{\right)}
\newcommand{\la}{\langle}
\newcommand{\ra}{\rangle}

\DeclareMathOperator{\spec}{{{Spec}}}

\newcommand{\st}{{\operatorname{st}}}
\newcommand{\ev}{{\operatorname{ev}}}
\newcommand{\ft}{\operatorname{ft}}
\newcommand{\vir}{\operatorname{vir}}

\newcommand{\pres}{\operatorname{pres}}

\newcommand{\vso}{\mathcal{O}^{\operatorname{vir}}}
\newcommand{\dertensor}{\otimes^{\text{der}}}
\newcommand{\on}{\operatorname}

\begin{document}
\title{Quantum $K$-Theory I: Foundations}
\author{Y.-P. Lee}
\address{Princeton University\\
        Department of Mathematics\\
    Fine Hall, Washington Road\\
        Princeton, NJ 08544-1000}
\email{yplee@math.princeton.edu}
\address{Department of Mathematics \\
        University of Utah \\
        Salt Lake City, Utah 84112-0090}
\email{yplee@math.utah.edu}
\thanks{Research partially supported by NSF grant DMS-0072547.}


\begin{abstract}
This work is devoted to the study of the foundations of \emph{quantum
$K$-theory}, a $K$-theoretic version of quantum cohomology theory.
In particular, it gives a deformation of the ordinary $K$-ring $K(X)$ of
a smooth projective variety $X$, analogous to the relation between quantum
cohomology and ordinary cohomology. This new quantum product also gives
a new class of Frobenius manifolds.
\end{abstract}

\maketitle





\section{Introduction} \label{S:0}
This work is devoted to the study of \emph{quantum $K$-theory},
a $K$-theoretic version of quantum cohomology theory. In particular, it
gives a deformation of the ordinary $K$-ring $K(X)$ of a smooth variety $X$,
analogous to the relation between quantum cohomology and ordinary cohomology.

In order to understand quantum $K$-theory, it is helpful to give a
brief review of quantum cohomology theory. Quantum cohomology
studies the intersection theory (of tautological classes) on
$\mbar_{g,n}(X,\beta)$, the moduli spaces of stable maps from
curves $C$ to a smooth projective variety $X$. The intersection
numbers, or Gromov--Witten invariants, look like
\[
  \int_{[\mbar_{g,n}(X,\beta)]^{\vir}} \st^*(\alpha)
    \prod_{i=1}^n \ev_i^*(\gamma_i) \psi_i^{d_i}
\]
where $\gamma_i \in H^*(X), \alpha \in H^*(\mbar_{g,n})$ and $\psi_i$ are the
cotangent classes. These notations will be defined in the next section.
For now these numbers are pairings between natural cohomology
classes and the (virtual) fundamental class of $\mbar_{g,n}(X,\beta)$.
These invariants at $genus(C)=0$ determines a deformation of the product
structure of the ordinary cohomology ring. The associativity of the quantum
cohomology ring structure is equivalent to the \emph{WDVV equation},
which is basically a degeneration argument
by considering a flat family of solutions parameterized by $\pp^1$ and
equating two different degeneration points of the solutions.

We propose to apply the same philosophy to $K$-theory.
The functorial interpretation of the integration over the (virtual) fundamental
classes $[\mbar_{g,n}(X,\beta)]^{\vir}$ is the push-forward of cohomologies
from $\mbar_{g,n}(X,\beta)$ to the point $\spec \cc$ (via the orientation
defined by virtual fundamental class). In fact the push-forward operation
can be performed in the relative setting for any proper morphisms.
With this in mind, we define \emph{quantum $K$-invariants} as
the $K$-theoretic push-forward to $\spec \cc$ of some natural vector
bundles on $\mbar_{g,n}(X,\beta)$ (via the orientation defined by the
 virtual structure sheaf). More precisely,
\[
  \chi \Bigl(\mbar_{g,n}(X,\beta), \vso_{\mbar_{g,n}(X,\beta)} \st^*(\alpha)
    \prod_{i=1}^n \ev_i^*(\gamma_i) \cl_i^{d_i} \Bigr)
\]
where $\gamma_i \in K(X), \alpha \in K(\mbar_{g,n})$ and $\cl_i$ are the
cotangent line bundles. The $K$-theoretic push-forward $\chi$
will be defined in Section~2.1 and Section~4.1 for algebraic and topological
$K$-theory respectively.
The quantum $K$-product is defined by quantum $K$-invariants and
a \emph{deformed} metric (equation \eqref{e:qm}).
The associativity of quantum $K$-product will be established by a
sheaf-theoretic version of WDVV-type argument. A considerable difference
in $K$-theory and intersection theory arises here.

Our main motivations for studying this theory come from two sources. The
first is to use this theory to study the geometry of the moduli space of maps.
For example, in joint works with R.~Pandharipande and with I.~Ciocan-Fontanine,
a computation of quantum $K$-invariants of $X=\pp^1$ are used to show that the
Gromov--Witten loci $\cap_i \ev^{-1}_i(pt_i) \subset \mbar_{0,n}(\pp^1, d)$
are not rational when $n$ is large. This thus answers a question raised in
\cite{LP} Section 2.2.

The other main motivation comes from the relation between Gromov--Witten
theory and integrable systems. There are many celebrated examples of this sort
in quantum cohomology theory. Two of them are most relevant.
The first one is Witten--Kontsevich's theory of two dimensional gravity
\cite{MK} \cite{EW}. This theory builds a link between (full)
Gromov--Witten theory of a point to KdV hierarchy. The second example is
Givental--Kim's theory: It states that genus zero Gromov--Witten theory
on flag manifolds is governed by quantum Toda lattices \cite{GK} \cite{BK}.
Our goal is therefore to extend this relation to quantum $K$-theory.
The analogous relation between (full)
quantum $K$-theory of a point and \emph{discrete} KdV hierarchy is yet to be
spelled out. See \cite{YL1} \cite{YL2} for some computations in this direction.
The relation between genus zero quantum $K$-theory on flag manifolds and
\emph{finite difference} Toda lattices turns out to be a very interesting one.
The construction of finite difference Toda lattices involved quantum
groups. The interested reader is referred to \cite{GL} for the details.

Besides these, we believe that the quantum $K$-invariants for symplectic
manifolds are symplectic invariants. This assertion, if verified, could have
some implications in symplectic geometry.
Also, this theory could be understood as one
manifestation of J.~Morava's quantum generalized cohomology theory \cite{JM}.

\begin{acknowledgements}
I am grateful to A.~Givental and R.~Pandharipande for much advice
in the course of this work. The referees' reports have greatly
improved the exposition of this work. In particular, the proof of
Proposition~\ref{p:contractions} is drastically simplified. Thanks
are also due to T.~Graber, G.~Kennedy, J.~Morava, Y.~Ruan and
V.~Srinivas for useful discussions. Part of the work is written
during my stay in NCTS, whose hospitality is appreciated.
\end{acknowledgements}

\section{Virtual structure sheaf}
The main purpose of this section is to define an element $\so^{\vir}$ in
the Grothendieck group of coherent sheaves on the moduli space of
stable maps to a smooth projective variety $X$.

\subsection{Preliminaries on algebraic $K$-theory}
The basic reference for algebraic $K$-theory discussed here is
\cite{sga6}.


Let $K^{\circ}(V)$ denote the Grothendieck group of locally free sheaves
on $V$ and $K_{\circ}(V)$ denote the Grothendieck group of coherent
sheaves on $V$. There is a cap-product
\[
  K^{\circ}(V) \otimes K_{\circ}(V) \to K_{\circ}(V)
\]
defined by tensor product
\[ [E]\otimes [F] \mapsto [E \otimes_{\so_V} F] , \]
making $K_{\circ}(V)$ a $K^{\circ}(V)$-module.

For any proper morphism $f: V' \to V$, there is a \emph{push-forward}
homomorphism
\[
  f_*: K_{\circ}(V') \to K_{\circ}(V)
\]
defined by
\[
 f_*([F]):= \sum_i (-1)^i [R^i f_* F].
\]

\begin{convention} If $f$ is an embedding, we will sometimes denote
$f_*([F])$ simply by $[F]$.
\end{convention}

A morphism $f:V' \to V$ is a \emph{perfect morphism} if there is an $N$
such that
\[
  Tor^V_i(\so_{V'}, F)=0
\]
for all $i>N$ and for all coherent sheaves $F$ on $V$.
For a perfect morphism $f:V' \to V$, one can define the \emph{pull-back}
\[
 f^*: K_{\circ}(V) \to K_{\circ}(V')
\]
by
\[ f^*([F])=\sum_{i=0}^d (-1)^i [Tor^{V}_i(\so_{V'}, F)]. \]
It is easy to see that regular embeddings and flat morphisms are perfect.

Let
\begin{equation} \label{e:gysin}
 \begin{CD}
   V' @>u>>V \\
  @VaVV @VbVV\\
  B' @>{v}>> B
 \end{CD}
\end{equation}
be a fibre square with $v$ a regular embedding. The \emph{refined Gysin map}
\[
 v^!: K_{\circ}(V) \to K_{\circ}(V')
\]
can be defined as follows. Since $v$ is a regular embedding,
$v_*(\so_{B'})$ has a finite resolution of vector bundles on $B$.
That is, $\so_{B'}$ determines an element $K^{\circ}_{B'}(B)$,
the Grothendieck group of vector bundles on $B$ whose homologies are
supported on $B'$. This then determines an element, denoted by
$[\so_{B'}^V]$, in $K^{\circ}_{V'}(V)$ by pull-back. Define
$K_{\circ}^{V'}(V)$ to be the Grothendieck group of the coherent sheaves
on $V$ with homologies supported on $V'$ (i.e.~exact off $V'$).
There is a canonical
isomorphism $\varphi: K_{\circ}^{V'}(V)\overset{\cong}{\to}
K_{\circ}(V')$. One can introduce a generalized cap product
\[
  K^{\circ}_{V'}(V) \otimes K_{\circ}(V) \overset{\otimes}{\to}
  K_{\circ}^{V'}(V) \overset{\varphi}{\to}  K_{\circ}(V')
\]
by composing the isomorphism $\varphi$ and the tensor product.
The refined Gysin map is:
\[
  v^!([F]) =  \varphi( [\so_{B'}^V \otimes_{\so_V} F]).
\]

\begin{remark} \label{r:1}
One can also apply the deformation to the normal cone argument to define
the Gysin pull-back. That is, one can replace the fibre square \eqref{e:gysin}
with
\[
 \begin{CD}
   V' @>u_C>> C_{V'} V \\
  @VaVV @VbVV\\
  B' @>{v_N}>> N_{B'} B
 \end{CD}
\]
where $C_{V'} V, N_{B'} B$ are the normal cone and normal bundle of
the embeddings $u, v$. The specialization to the normal cone argument
works for $K$-theory as well, i.e.~one has a map
\[
  \sigma: K_{\circ}(V) \to K_{\circ}(C_{V'} V)
\]
such that the Gysin map can be defined by combining $\sigma_V$
with the Gysin map $v_N^!$ of the above fibre square. That is
\[
  v^!=v_N^! \circ \sigma.
\]
\end{remark}

\begin{remark} \label{r:2}
If $b$ is also a regular embedding, then $v^! [\so_V]= b^!
[\so_{B'}]$, which will also be denoted $[\so_{B'}]
\dertensor_{\so_V} [\so_V']$ to emphasize the symmetry.
\end{remark}

If $v: B' \to B$ can be factored through a regular embedding $i:B'
\to Y$ followed by a smooth morphism $p: Y \to B$, i.e.~$v$ is a
\emph{local complete intersection morphism}. One can form the
following fibre diagram
\[
 \begin{CD}
  V' @>j>> Y' @>q>> V \\
  @VVV   @VVV   @VVV \\
  B' @>i>>  Y @>p>> B .
 \end{CD}
\]
Then $q$ is smooth and we define $v^!([F]):=i^! (q^* [F])$. This
definition is actually independent of the factorization.

\begin{convention}
To simplify the notation, we will often denote the class $[F]$ simply by $F$
and use the operations (like $\otimes$) as if $[F]$ is a complex of sheaves.
\end{convention}

\subsection{Moduli spaces of stable maps}
This section serves only the purpose of fixing the notations. The reader
should consult \cite{FP} \cite{BM} for details.

\begin{definition*}
A \emph{stable map} $(\mathcal{C};x,f)$ from a
pre-stable curve $(C;x_1, \ldots, x_n)$ of genus $g$ curves
with $n$ marked points to
$X$, which represents a class $\beta \in H_2(X)$, is a morphism
$f: C \to X$ satisfying:

1. The homological push-forward of $C$ satisfies ${f}_*([C])=\beta$.

2. \emph{Stability}: The following
inequality holds on the normalization $E'$ of each irreducible component $E$:
\[ \text{$f(E)$ is a point}
\Rightarrow 2g(E')+\#\{\text{special points on $E'$}\} \ge 3
\]
\end{definition*}

It is obvious from the definition that the notion of stable maps is a natural
generalization of the notion of stable curves. The essential point being that
they only allow only \emph{finite} automorphisms. This makes
the following theorem possible:

\begin{theorem} {\rm (\cite{MK})}
The stack $\mbar_{g,n}(X,\beta)$ of stable maps to $X$, a smooth projective
scheme of finite type over $\cc$, is a proper Deligne-Mumford stack.
\end{theorem}

Between the various spaces in Gromov--Witten theory there are many interesting
morphisms. The first such class of morphisms is called the \emph{forgetful}
morphisms
\[ \ft_k : \mbar_{g,n+1}(X,\beta) \to \mbar_{g,n}(X,\beta) \]
by forgetting $k$-th marked points and stabilizing if necessary.

The second class is called the \emph{stabilization} morphisms
\[ \st:\mbar_{g,n}(X,\beta) \to \mbar_{g,n}, \]
where $\mbar_{g,n}$ is the moduli
stack of stable curves. It is defined by forgetting
the map and retaining the pre-stable curve, then stabilizing the curve.
Of course one can just forget the map without stabilization.
This leads to the \emph{prestabilization} morphism \cite{BM}
\[ \pres: \mbar_{g,n}(X,\beta) \to \mathfrak{M}_{g,n},  \]
with $\mathfrak{M}_{g,n}$ being the moduli stack of prestable curves.

The third class is called the \emph{evaluation} morphisms
\[   \ev_i :\mbar_{g,n}(X,\beta) \to X \]
by evaluating the stable map at the $i$-th marked point.
That is, $\ev_i(f)=f(x_i)$.

Another object which we will use quite often is the \emph{universal cotangent
line bundle} $\cl_i$ on $\mbar_{g,n}(X,\beta)$. It is defined as
$\cl_i=x_i^*(\omega_{\mathcal{C}/M})$, where
$\omega_{\mathcal{C}/M}$ is the relative dualizing sheaf of the universal
curve $\mathcal{C} \to \mbar_{g,n}(X,\beta)$ and $x_i$ are the marked points.

There are combinatorial generalizations of $\mbar_{g,n}(X,\beta)$
by incorporating the modular graphs $\tau$ with degree $\beta$
\cite{BM}. The idea is to use the graphs to keep track of the
combinatorics of degeneration of curves. This bookkeeping
simplifies the notations in the proof of functorial properties of
the virtual structure sheaf. We briefly recall the definitions
here.

A \emph{graph} $\tau$ is a quadruple $(F_{\tau},V_{\tau},j_{\tau},\p_{\tau})$
where $F_{\tau}$ (\emph{flags}) and $V_{\tau}$ (\emph{vertices}) are finite
sets, $j : F_{\tau} \to F_{\tau}$ is an involution and
$\p : F_{\tau} \to V_{\tau}$ is a map.
\[ S_{\tau}:= \{f \in F_{\tau}| jf=f\} \]
is the set of \emph{tails} (or \emph{marked points}) and
\[ E_{\tau}:= \{ \{f_1,f_2\} \subset F_{\tau} | f_2 =j f_1, f_2 \ne f_1 \} \]
is the set of \emph{edges}.
A \emph{modular graph} is a pair $(\tau,g)$, where $\tau$ is
a graph and $g: V_{\tau} \to \zz_{\ge 0}$ is a map. We call $g(v)$ the
\emph{genus} of the vertex $v$. A modular graph is \emph{stable} if
\[
  2 g(v) + \# F_{\tau}(v) \ge 3
\]
for all vertices $v$, where $F_{\tau}(v):= j_{\tau}^{-1}(v)$ and $\# F(v)$ is
the valence of $v$.

One defines the moduli stack of prestable curves associated to a modular graph
$(\tau,g)$ as
\[
  \mathfrak{M}_{\tau,g} :=
  \prod_{v \in V_{\tau}} \mathfrak{M}_{g(v), \# F_{\tau} (v)}
\]
where $\mathfrak{M}_{g,n}$ is the usual moduli stack of prestable curves.
The moduli stack of stable curves $\mbar_{\tau,g}$ associated to a
\emph{stable} modular graph $(\tau,g)$ is defined similarly
\[
  \mbar_{\tau,g} :=
  \prod_{v \in V_{\tau}} \mbar_{g(v), \# F_{\tau} (v)}.
\]
Geometrically, the moduli stack of (pre)stable curves associated to a
(connected) modular graph $(\tau,g)$ which has at least one edge consists of
(connected) singular curves with singularity prescribed by $\tau$.
In other words, the stack of universal curves over the stack of the (pre)stable
curves has singular generic fibres. For details, see \cite{BM}.

Let $(\tau,g)$ be a modular graph and $\beta: V_{\tau} \to H_2(X)$ be the
\emph{degree}.
\footnote{We have abused the notations by using $\beta$ (resp. $g$)
sometimes as map $\beta: V_{\tau} \to H_2(X)$ ($g: V_{\tau} \to \zz_{\ge 0}$)
and sometimes as the total degree $\sum_v \beta(v)$ (total genus).
It should be clear in the context which is which.}
$(\tau,g,\beta)$ is a \emph{stable modular graph with degree} if
\[
  \beta(v)=0 \Rightarrow 2 g(v) + \# F(v) \ge 3.
\]
The moduli stack of stable maps $\mbar(X,\tau,g,\beta)$
\footnote{For brevity, $\tau$ will sometimes stand for $(\tau,g)$
or $(\tau,g, \beta)$.}
associated to $(\tau,g,\beta)$ is defined
by the following three conditions:
\begin{enumerate}
\item If $\tau$ is a graph with only one vertex $v$, and set of flags $F_v$
  are all tails, then
 \[ \mbar(X,\tau,g,\beta):= \mbar_{g(v),\# F_v}(X,\beta).
 \]
\item
 \[
  \mbar(X,\tau_1 \times \tau_2) := \mbar(X,\tau_1) \times \mbar(X,\tau_2),
 \]
where $\tau_1 \times \tau_2$ is the disjoint union of two graphs and
$\tau$ stands for $(\tau, g,\beta)$, by abusing the notation.
\item If $\sigma$ is obtained from $\tau$ by cutting an edge, then
$\mbar(X,\tau)$ is defined by the following fibre square
\begin{equation} \label{e:1}
 \begin{CD}
  \mbar(X,\tau) @>>> \mbar(X,\sigma) \\
  @VVV      @VV{\ev_i \times \ev_j}V \\
  X @>\Delta>> X\times X,
 \end{CD}
\end{equation}
where $\Delta$ is the diagonal morphism, and $i,j$ are the corresponding
marked points where the edge is cut.
\end{enumerate}

It is easy to see that the morphisms and operations on $\mbar_{g,n}(X,\beta)$
generalize to $\mbar(X,\tau)$. See \cite{BM} for the details.

\subsection{The construction of virtual structure sheaf} \label{S:2}
We shall now introduce the notion of \emph{virtual structure sheaf},
whose construction is parallel to that of virtual fundamental
class as defined in \cite{LT} \cite{BF}.

Let $V \to B$ be a morphism from a DM-stack to a smooth Artin stack
with constant dimension.
A \emph{relative perfect obstruction theory} is a homomorphism $\phi$ in the
derived category (of quasi-coherent sheaves bounded from above) of a two term
complex of vector bundles $E^{\bullet} = [E^{-1} \to E^0]$ to the relative
cotangent complex $L^{\bullet}_{V/B}$ of $V \to B$ such that:

1. $H^0(\phi)$ is an isomorphism.

2. $H^{-1}(\phi)$ is surjective.

If $V$ admits a global embedding $i: V \to Y$ over $B$ into a smooth
Deligne-Mumford stack, then the two term cutoff of $L^{\bullet}_{V/B}$
at $-1$ is given by the first two terms in second exact sequence of
relative K\"ahler differential: $[I/I^2 \to  i^*(\Omega_{Y/B})]$, where $I$
is the relative ideal sheaf of $V$ in $Y$ over $B$.

Let us now recall the definition of \emph{relative intrinsic
normal cone}. Choose a closed $B$-embedding of $V$, $i: V
\hookrightarrow Y$, to a smooth DM stack $Y$ over $B$. Since the
normal cone $C$ is a $T_Y$-cone over $B$ (\cite{BF}), one could
associate a cone stack $\mathfrak{C}_{V|B}:=[C/i^*T_Y]$, which is
independent of the choice of embedding. $E^{\bullet} \to
L^{\bullet}_V/B$ being a perfect obstruction theory implies
(\cite{BF} Theorem~4.5) that $\mathfrak{C} \to [E_1/E_0]$ is a
closed embedding, which induces a cone $C_1 \subset E_1$. Here
$E_0 \to E_1$ is the dual complex of $E^{-1} \to E^0$. Note that
even if $V$ does not allow a global embedding, one can still
define the intrinsic normal cone by (\'etale) local embedding. It
is proved in \cite{BF} that the local construction glues together
to form a global cone stack.

\begin{definition*}
Let $E^{\bullet} \to L^{\bullet}_{V/B}$ be a perfect obstruction theory.
The virtual structure sheaf $\so^{\vir}_{[V,E]}$ for a (relative) perfect
obstruction theory is an element in $K_{\circ}(V)$ defined as the
pull-back of $\so_{C_1} \in K_{\circ}(E_1)$ by the zero section
$0: V \to E_1$,
i.e.,
\[
 \so^{\vir}_{[V,E]}= \sum_i (-1)^i [Tor_i^{E_1}(\so_{V} \otimes \so_{C_1})].
\]
\end{definition*}

\begin{remark}
One could also consider this as an element in
the derived category of $V$ defined to be
\[ \so^{\vir}_{V} := \so_{C_1} \dertensor_{\so_{E_1}} \so_{V},\]
where $\dertensor$ means the derived tensor.
Either way, virtual structure sheaf is in general not a bona fide sheaf.
\end{remark}

Apply the above setting to Gromov--Witten theory.
Let $\mathfrak{M}_{\tau}$ be the corresponding moduli
stack associated with modular graph $\tau$.
Recall that there is a prestabilization morphism
$\pres: \mbar(X,\tau) \to \mathfrak{M}_{\tau}$, which will be our
structure morphism $V \to B$.
Consider the universal stable map
\begin{equation} \label{e:usm}
 \begin{CD}
  \mathcal{C} @>f>> X \\
  @V{\pi}VV \\
  \mbar(X,\tau).
 \end{CD}
\end{equation}
One then get a complex $R^{\bullet} \pi_* f^* T_{X}$,
which can be realized as a two term complex $[E_0 \to E_1]$
($E_i^{\vee} = E^{-i}$) of vector bundles on $\mbar_{g,n}(X,\tau)$
(\cite{KB}). In fact, there is a homomorphism
\begin{equation*}
 \phi: (R^{\bullet} \pi_* f^* T_{X})^{\vee} \to
 L^{\bullet}_{\mbar(X,\tau)/\mathfrak{M}_{\tau}}
\end{equation*}
which is a relative perfect obstruction theory.
One therefore defines the virtual structure sheaf $\so_{\mbar(X,\tau)}^{\vir}$.
Note that in the present case $\mbar(X, \tau)$ admits a global embedding
(Appendix A of \cite{GP}). One might identify
$L^{\bullet}_{\mbar(X,\tau)/\mathfrak{M}_{\tau}}$ with its two-term cutoff.


\subsection{Basic properties of virtual structure sheaf}
Some basic properties of virtual structure sheaf is established here as a
$K$-theoretic version of \cite{BF}.
Because $B$ will always be $\mathfrak{M}_{\tau}$ in all our applications,
the readers should feel free to make this assumption, although all the results
will hold for $B$ any smooth equidimensional Artin stack.

The first step to check the consistency of this definition is the following
two propositions:

\begin{proposition} \textbf{(Consistency)} \label{p:consistency}
The virtual structure sheaf is independent of the choice of the global
resolutions $E^0 \to E^1$. That is, if $E^0 \to E^1$ is quasi-isomorphic to
$F^0 \to F^1$, then $\vso_{[V,E]}=\vso_{[V,F]}$.
\end{proposition}

\begin{proof}
The proof of this proposition is exactly the same as the proof of
\cite{KB} Proposition~5.3.
\end{proof}

\begin{proposition}  \textbf{(Expected dimension)} \label{p:ed}
If $V$ has the \emph{expected dimension} $\dim B+
\on{rank} (E^{\bullet})$, then $V$ is a local complete intersection and
the virtual structure sheaf is equal to ordinary structure sheaf. Note that
the dimension is called ``expected dimension'' because it is the dimension
in the (``generic'') case of no obstruction.
\end{proposition}

\begin{proof}
The question can be reduced to the absolute case by the argument given in
\cite{GP} Appendix B. Thus we will set $B=pt$.
Since the question is local (in \'etale topology), one may assume that
\[ X=\spec R, \quad R=\cc[x_1,x_2,\ldots,x_r]/(g_1,\ldots,g_s). \]
In fact, it is enough to consider the local ring
\[
  R_p=\left(\dfrac{\cc[x_1,x_2,\ldots,x_r]}
  {(g_1,\ldots,g_s)} \right)_{(x_1,\ldots,x_r)}
\]
with maximal ideal $\mathfrak{m}$.
Denote $I=(g_1,\ldots,g_s)$, then the two-term cut-off of cotangent
complex is
\[ 
 I/I^2 \to \Omega_{\cc[x_1,\ldots,x_r]} \otimes_{\cc[x_1,\ldots,x_r]} R . 
\]
This induces a complex of $R_p$ modules
\begin{equation*}
  I/\mathfrak{m}I \overset{\varphi}{\to} \mathfrak{m}/ \mathfrak{m}^2.
  \tag{*}
\end{equation*}
One may assume that $g_i$ contains no linear terms in $x_i$, or otherwise
the presentation of $R_p$ could be changed to eliminate terms linear in $x_i$.
That is, $I \subset \mathfrak{m}^2$ and the above $\varphi$ is a zero map.
This implies that the homology of the above complex is
\[ 
  H^1(*) =I/\mathfrak{m}I, \quad H^0(*)=\mathfrak{m}/\mathfrak{m}^2.
\]
Because $E^{\bullet} \to L^{\bullet}$ is a perfect obstruction theory,
\[ 
  H^0(*)= H^0(E^{\bullet \vee} \otimes \cc), \quad
  H^1(*) \hookrightarrow H^1(E^{\bullet \vee} \otimes \cc).
\]
Therefore
\[
 \on{rank}(E^{\bullet})= h^0(E^{\bullet}\otimes \cc) -
  h^1(E^{\bullet}\otimes \cc) \le h^0(*) - h^1(*) = r-s
\]
where the last equality is obvious. This shows, combining with
the assumption, that $\dim \spec R_p \le r-s$. However, $\dim \spec R_p$ is
obviously greater or equal to $r-s$. Thus $\dim \spec R_p = r-s$ and
$\{ g_i \}$ form a regular sequence. Therefore
$V$ is a local complete intersection.

On the other hand, $\dim \spec R_p = r-s$ implies that
$h^1(E^{\bullet \vee} \otimes \cc) = s$ and $E^{\bullet} \to L^{\bullet}$
is an isomorphism. Hence $C=E_1$,
and virtual structure sheaf is the ordinary structure sheaf on $V$.
\end{proof}

The following two propositions form the technical heart of the virtual
structure sheaves.

Consider a fibre square
\begin{equation} \label{e:3}
 \begin{CD}
  V' @>u>>  V \\
  @VaVV   @VbVV \\
  B' @>v>>  B
 \end{CD}
\end{equation}
where $V', V$ are separated DM stacks
and $B',B$ are smooth (Artin) stacks of constant dimensions.

\begin{proposition} \label{p:pullback} \textbf{(Pull-back)}
Suppose that $E^{\bullet} \to L^{\bullet}_{V/B}$ is a perfect
obstruction theory for $V \to B$, and $v$ is flat or is a regular embedding.
Then $u^*E^{\bullet} \to L^{\bullet}_{V'/B'}$ is a perfect obstruction theory
for $V' \to B'$ and
\[
  \vso_{[V',u^*E^{\bullet}]} = v^! \vso_{[V,E^{\bullet}]}.
\]
\end{proposition}

Let $E$ be a perfect relative obstruction theory for $U$ over $B$
and $E'$ be a perfect relative obstruction theory for $U'$ over $B$, where
$U, U'$ are DM stacks and $B$ is a smooth (Artin) stack of constant dimension.
\footnote{The superscript $\bullet$ in $E^{\bullet}$ will be omitted for
simplicity.}
Consider the following fibre square of DM $B$-stacks
\begin{equation} \label{e:4}
 \begin{CD}
  U' @>{\mu}>> U \\
  @V \alpha VV @V \beta VV \\
  V'@>{\nu}>> V.
 \end{CD}
\end{equation}
Suppose that the morphism $\nu$ is a local complete intersection morphism of
$B$-stacks that have finite unramified diagonal over $B$. Recall
(\cite{BF} Sect.~7) that $E$ and $E'$ are said to be \emph{compatible
over $\nu$} if there exists a homomorphism of distinguished triangles
\[
\begin{CD}
 \mu^* E @>>> E' @>>> \alpha^* L_{V'/V} @>>> \mu^*E[1] \\
 @VVV   @VVV    @VVV        @VVV \\
 \mu^* L_{U/B} @>>> L_{U'/B} @>>> L_{U'/U} @>>> \mu^* L_{U/B}[1].
\end{CD}
\]
The compatibility of obstruction theory means that we have the following
short exact sequence of vector bundle stacks
\begin{equation} \label{e:compatibility}
 \alpha^* \mathfrak{N}_{V'/V} \to [E'_1/E'_0] \to [\mu^* E_1/\mu^* E_0],
\end{equation}
where $\mathfrak{N}_{V'/V} := h^1/h^0 (T^{\bullet}_{V'/V})$. (In general
a vector bundle stack can be locally represented as a quotient of
a vector bundle by another vector bundle.) In particular
when $\nu$ is a local complete intersection morphism,
one can find a smooth stack $Y$ and $V' \to V$ factors through $Y$:
\[
  V' \overset{i}{\to} Y \overset{p}{\to} V
\]
such that $i$ is regular embedding and $p$ is smooth. Then
\[
  \mathfrak{N}_{V'/V} := [N_{V'/Y} / i^* T_{Y/V}].
\]

\begin{proposition} \label{p:functoriality} \textbf{(Functoriality)}
If $E'$ and $E$ are compatible over $\nu$, then
\[
 \nu^! \vso_{[U,E]} = \vso_{[U',E']}
\]
when $\nu$ is smooth or $V'$ and $V$ are smooth.
\end{proposition}

The rest of this subsection is devoted to the proof of the above two
propositions. It can be skipped without losing the logic flow.

Consider the following fibre diagram
\begin{equation} \label{e:5}
 \begin{CD}
 C_{X_1}Y \times_Y C_{X_2}Y @>>> C_{X_2}Y \times_Y X_1 @>>> C_{X_2}Y \\
 @VVV   @VVV    @VVV \\
 C_{X_1}Y \times_Y X_2 @>>> Z @>>> X_2 \\
 @VVV   @VVV    @ViVV \\
 C_{X_1}Y @>\rho>> X_1 @>j>> Y.
\end{CD}
\end{equation}
where $i,j$ are local embeddings.

\begin{lemma} \label{l:1}
$[\so_{C_1}]=[\so_{C_2}]$ in $K_{\circ}(C_{X_1}Y \times_Y C_{X_2}Y)$, where
$C_1 := C_{C_{X_1}Y \times X_2} (C_{X_1}Y)$ and $C_2$ is similarly defined
by switching indices $1$ and $2$. $C_i$ are naturally embedded in
$C_{X_1}Y \times_Y C_{X_2}Y$.
\end{lemma}

\begin{proof}
The proof of Proposition~4 in \cite{AK} also proves this lemma:
It is constructed there two principal Cartier divisors $D$ and $E$ in
$M^{\circ}_{X_1} Y \times_Y M^{\circ}_{X_2} Y$ (with notations like
\cite{WF} Ch.~5). The proof there shows that $i_E^{!} [\so_D] = [\so_{C_1}]$
and $i_D^{!} [\so_E] = [\so_{C_2}]$. By Remark~\ref{r:2},
$[\so_{C_1}]=[\so_{C_2}]$.
\end{proof}

Now apply the above to the diagram \eqref{e:4}. (Similar arguments applies to 
diagram \eqref{e:3}.) Assuming that $\nu$ is a regular embedding and $\beta$ 
is an embedding, one has
\[
 \begin{CD}
 N \times_V D @>>> D' @>>> D \\
 @VVV   @VVV    @VVV \\
 N \times_V U @>>> U' @>\mu >> U \\
 @VVV   @V \alpha VV    @V \beta VV \\
 N @>\rho>> V' @> \nu >> V,
\end{CD}
\]
where $D=C_U V$, $D'=C_{U'}V'$ and $N=N_{V'}V$. The Lemma~\ref{l:1} states that
\[
  [\so_{C_{D'}D}] = [C_{N \times_V U} N] = [\rho^* \so_{D'}],
\]
where the $\rho^*$ is the pull-back on cones induced from $\rho$.
Let $0: D' \to N \times_V D$ be the zero section.
By the deformation to normal cone argument
\[
  \nu^! [\so_D] = 0^* [\so_{C_{D'}D}] \in K_{\circ}(D').
\]
Also
\[
  0^* [\rho^* \so_{C_{U'}V'}] = [\so_{C_{U'}V'}]  \in K_{\circ}(D').
\]
This implies that

\begin{lemma} \label{l:2}
\[
 \nu^! [\so_D] = [\so_{C_{U'}V'}].
\]
\end{lemma}

For the readers who are familiar with the arguments in \cite{BF},
the above lemmas easily implies the propositions. In the following we
reproduce their arguments in $K$-theory.

Let us start with the proof of {Functoriality}. For simplicity, let $B=pt$.
We will follow closely the proof of Proposition~5.10 in \cite{BF}.
By Proposition~\ref{p:consistency} we may choose the global resolution
$[E'_0 \to E'_1]$ of $E'$ such that the right square of the following diagram
\[
 \begin{CD}
  0 @>>> F_0 @>>> E'_0 @>\phi_0>> \mu^* E_0 @>>> 0\\
  && @VVV   @VVV    @VVV    \\
  0 @>>> F_1 @>>> E'_1 @>\phi_1>> \mu^* E_1 @>>>0
 \end{CD}
\]
commutes and $\phi_i$ are surjective. $F_i:= \ker(\phi_i)$
are vector bundles. This induces short exact sequence
of vector bundle stacks
\[
  [F_1/F_0] \to [E'_1/E'_0] \to [\mu^* E_1/\mu^* E_0] .
\]
The compatibility of obstruction theory means that we have the following
short exact sequence of vector bundle stacks \eqref{e:compatibility}
\[
 \alpha^* \mathfrak{N}_{V'/V} \to [E'_1/E'_0] \to [\mu^* E_1/\mu^* E_0].
\]
This implies that $[F_1/F_0] \cong \alpha^* \mathfrak{N}_{V'/V}$.
Let $C_1=\mathfrak{C}_U \times_{\mathfrak{E}} E_1$
and $C'_1=\mathfrak{C}_U' \times_{\mathfrak{E}'} E'_1$. Then $\vso_{[U,E]}
=0_{E_1}^* \so_{C_1}$ and $\vso_{[U',E']}=0_{E'_1}^* \so_{C'_1}$.

When $\nu$ is smooth or a regular embedding, Proposition~3.14 of \cite{BF}
implies that we have a fibre square
\[
 \begin{CD}
  C'_1 @>>> \mu^* C_1 \\
  @VVV  @VVV \\
  E'_1 @>>> \mu^* E_1.
 \end{CD}
\]
If $\nu$ is smooth, $\nu^!=\mu^*$ and we have
\[
  \nu^! \vso_{[V,E]} = \nu^! 0^*_{E_1} [\so_{C_1}]
  = 0^!_{\mu^*E_1}[\mu^* \so_{C_1}]
  = 0^!_{E'_1} [\so_{C'_1}] =\vso_{[V',E']}.
\]
If $\nu$ is a regular embedding of smooth stacks, then
\[
 \nu^! \vso_{[V,E]} = \nu^! 0^*_{E_1} [\so_{C_1}] = 0^!_{E'_1} \nu^![\so_{C_1}]
  = 0^!_{E'_1} [\so_{C'_1}] =\vso_{[V',E']}
\]
by Lemma~\ref{l:2}. In general, when only smoothness of $V'$ and $V$ are
assumed, one factors $v$ as
\[
 \begin{CD}
  U' @>>> U' \times U @>>> U \\
  @VVV  @VVV    @VVV \\
  V' @>>> V' \times V @>>> V.
 \end{CD}
\]
Apply the smooth case to the right square and the regular embedding case
to the left square.

One can generalize the functoriality property
to the case $\beta$ an arbitrary morphism by applying the arguments
in \cite{BF} Lemma~5.9. That is, to replace
the morphism $j:U\to V$ by $U \to V \times Y$ such that $U \to Y$ (and
therefore $U \to V \times Y$) is an embedding. The fibre square \eqref{e:4}
will look like
\[
 \begin{CD}
  U' @> \mu >> U \\
  @VVV @VVV \\
  V'\times Y @> \nu >> V\times Y.
 \end{CD}
\]
This ends the proof of Proposition\ref{p:functoriality}.

The proof of Proposition~\ref{p:pullback} is similar.
Let $E^{\bullet}=E^{-1} \to E^0$. In the case $v$ is flat,
by the standard properties of (intrinsic) normal cones
$\mathfrak{C}_{V'/B'}= \mathfrak{C}_{V/B} \times_B B'$
(\cite{BF} Proposition~7.1) and $C'_1 \subset u^* E_1$ is the pull-back of
$C_1$. Therefore $\so_{C'_1}= u^* \so_{C_1} = v^! \so_{C_1}$.
Since $0^*$ commutes with $v^!$,
\[
 \vso_{[V',u^*E]}= 0^* \so_{C'_1}= 0^* v^! \so_{C_1} = v^! 0^* \so_{C_1}
  = v^! \vso_{[V,E]}.
\]
If $v$ is a regular embedding, apply Lemma~\ref{l:2} to get
$[\so_{C'_1}]= v^! [\so_{C_1}]$ in $K$-group and proceed as above.

\section{Axioms of virtual structure sheaves in quantum
$K$-theory} \label{s:3}
The axioms listed here are $K$-theoretic version of the Behrend--Manin
axioms \cite{BM}.

\subsection{Mapping to a point}
Suppose that $\beta=0$. In this case $\mbar(X,\tau) = \mbar_{\tau} \times X$
and 
\[
 R^1 \pi_* f^* T_X= R^1 \pi_* \so_{\mathcal{C}} \boxtimes T_X,
\]
where the notations are defined in \eqref{e:usm}.
Then
\begin{proposition}
\begin{equation*}
 \vso_{\mbar(X,\tau)} = \lambda_{-1}(R^1 \pi_*
  \so_{\mathcal{C}} \boxtimes T_X)^{\vee},
\end{equation*}
where $\lambda_{-1}(F):=\sum_i (-1)^i \wedge^i F$ is the alternating sum of
the exterior power of $F$.
\end{proposition}

\begin{proof}
The (dual of the) obstruction theory
$E_{\bullet}$ is quasi-isomorphic to $R^0 \pi_*  \so_C \boxtimes T_X \to
R^1 \pi_*  \so_C \boxtimes T_X$, which is a
complex of locally free sheaves with zero differential.
Since $\mbar_{\tau} \times X$ is smooth, 
$C_1= \on{Image} (E_0 \to E_1)= \mbar_{\tau} \times X$.
Now use the self-intersection formula \cite{FL} Chapter~VI.
\end{proof}

\subsection{Products}
\begin{proposition}
\[
   \vso_{\mbar(X,\tau_1) \times \mbar(X,\tau_2)}
 = \vso_{\mbar(X,\tau_1)} \boxtimes \vso_{\mbar(X,\tau_2)},
\]
here $\vso_{\mbar(X,\tau_1) \times \mbar(X,\tau_2)}$ makes sense as
\[
 \mbar(X,\tau_1) \times \mbar(X,\tau_2)= \mbar(X,\tau_1 \times \tau_2).
\]
See Section~2.2.
\end{proposition}

\begin{proof}
It is obvious that the perfect obstruction theory of the product is the
product of perfect obstruction theories on each component.
\end{proof}

\subsection{Cutting edges}
Let $\sigma$ be a modular graph obtained from $\tau$ by cutting an edge
as in the fibre diagram \eqref{e:1}. In this case $\mathfrak{M}_{\sigma}
=\mathfrak{M}_{\tau}=:\mathfrak{M}$. Consider the fibre square
\[
 \begin{CD}
  \mbar(X,\tau) @>>> \mbar(X,\sigma) \\
  @VVV      @VV{\ev_i \times \ev_j}V \\
  \mathfrak{M} \times X @>\Delta>> \mathfrak{M} \times X\times X,
 \end{CD}
\]
where $i,j$ are the markings created by the cut.

The cutting edges axiom says that
\begin{proposition}
\[
 \Delta^! \vso_{\mbar(X,\sigma)}   = \vso_{\mbar(X,\tau)}.
\]
\end{proposition}

\begin{proof} It is proved in \cite{KB} that the perfect obstruction
theories of ${\mbar(X,\sigma)}$ and ${\mbar(X,\tau)}$ are compatible with
respect to $\Delta$. It remains to apply Proposition~\ref{p:functoriality}.
\end{proof}

\subsection{Forgetting tails}
Let $\sigma$ be obtained from $\tau$ by forgetting a tail (marked point)
and let
\[ \pi: \mbar(X,\tau) \to \mbar(X,\sigma)
\]
be the forgetful map. Then
\begin{proposition}
\[
  \pi^* \vso_{\mbar(X,\sigma)} = \vso_{\mbar(X,\tau)}.
\]
\end{proposition}

\begin{proof}
The proof in \cite{KB} p.~610 shows that
\[
 \begin{CD}
  \mathfrak{C}_{\mbar(X,\tau)/\mathfrak{M}_{\tau}} @>\cong>>
   \pi^* \mathfrak{C}_{\mbar(X,\sigma)/\mathfrak{M}_{\sigma}} \\
  @VV{\cap}V        @VV{\cap}V \\
  \mathfrak{E}(X,\tau) @>\cong>> \pi^* \mathfrak{E}(X,\sigma)
 \end{CD}
\]
form a commutative diagram. Therefore,
\[
  \vso_{\mbar(X,\tau)} = \vso_{[\mbar(X,\tau), \pi^* E(X,\sigma)]}.
\]
Now use the pull-back property Proposition~\ref{p:pullback}
\[
  \pi^* \vso_{\mbar(X,\sigma)} = \vso_{[\mbar(X,\tau), \pi^* E(X,\sigma)]}.
\]
This completes the proof.
\end{proof}

\begin{remark}
The above propositions in this section are the first four of five axioms
on virtual fundamental classes in quantum cohomology
theory listed in \cite{BM} .
The last one, \emph{isogenies} axiom, is not literally true in
$K$-theory. The reason is that even though one can identify fundamental
classes of two cycles by a birational morphism, it is not true that the
structure sheaves of these two cycles are identified.
We therefore have to modify the isogenies axiom in quantum $K$-theory.
Four types of isogenies are stated in $K$-theory in three axioms.
The contractions axiom has to be modified,
or ``quantized'', in $K$-theory. (The meaning of ``quantization''
will be clear in the next section.)
\end{remark}

\subsection{Fundamental classes}
Let $\sigma$ be obtained from $\tau$ by forgetting a tail. The degree of
$\sigma$ is induced from $\tau$ in the obvious way. Consider the commutative
diagram
\[
 \begin{CD}
  \mbar(X,\tau) @>>> \mbar(X,\sigma) \\
  @VVV      @VVV    \\
  \mbar_{\tau}  @>{\Phi}>>  \mbar_{\sigma}
 \end{CD}
\]
which induces a morphism
\[
 \Psi: \mbar(X,\tau) \to \mbar_{\tau} \times_{\mbar_{\sigma}} \mbar(X,\sigma).
\]

\begin{proposition}
\[
  \Psi_*(\vso_{\mbar(X,\tau)}) = \Phi^! \vso_{\mbar(X,\sigma)}
\]
\end{proposition}

\begin{proof}
The proof is the same as the proof in \cite{KB} pp.~611-3, with the
(bivariant) cycle classes changed to (bivariant) $K$-classes. We summarize
the spaces involved in the following diagram:
\[
 \xymatrix{\mbar(X, \tau) \ar[rr]^{\Psi} \ar[d]
  &&\mbar_{\tau} \times_{\mbar_{\sigma}} \mbar(X,\sigma) \ar[r] \ar[d]
  &\mbar(X,\sigma) \ar[d]^a \\
  \mathcal{C} \ar[rr]^l \ar[rd] \ar@/^2pc/[rrr]^{\pi'}
  &&\mbar_{\tau} \times_{\mbar{\sigma}} \mathfrak{M}_{\sigma} \ar[r]^m \ar[dd]
  &\mathfrak{M}_{\sigma} \ar[dd]^{s_{\sigma}} \\
  &\mathfrak{M}_{\tau} \ar[rd]^{s_{\tau}} &&\\
  &&\mbar_{\tau} \ar[r]^{\Phi} &\mbar_{\sigma}}
\]
where the squares are cartesian and
$\pi': \mathcal{C} \to \mathfrak{M}_{\sigma}$ is the universal curve.
Since $\pi'$, $m$ and $\Phi$ are representable, proper and flat, they have
natural orientations $[\pi'], [m], [\Phi]$ (which are bivariant $K$-classes
\cite{FM}). Since $s_{\sigma}$ is flat, $[m] = s_{\sigma}^*[\Phi]$.
Furthermore, since $l$ is just blow-down of some rational curves,
$l_* \so_{\mathcal{C}} =
\so_{\mbar_{\tau} \times_{\mbar{\sigma}} \mathfrak{M}_{\sigma}}$.
\footnote{Alternatively, one can argue this in the following way. It is
known that if $f:U \to V$ is a birational surjective morphism between
smooth schemes, then $R^i f_* \so_U = 0$ for $i \ge 1$ and
$R^0 f_* \so_U =\so_V$. This property is preserved under flat base change,
therefore it holds for stacks as well. That implies $l_* \so = \so$.}
Therefore $l_*[\pi']=s_{\sigma}^*[\Phi]$. Then
\begin{align*}
 \Phi^![\vso_{\mbar(X,\sigma)}]
 =&a^* s_{\sigma}^*[\Phi] \vso_{\mbar(X,\sigma)} \\
 =&a^* l_*[\pi'] \vso_{\mbar(X,\sigma)} \\
 =&\Psi_* {\pi'}^! \vso_{\mbar(X,\sigma)} \\
 =&\Psi_* \vso_{\mbar(X,\tau)}
\end{align*}
where the last equality follows from Proposition~\ref{p:functoriality}.
\end{proof}

\subsection{Isomorphisms} Suppose that $\sigma$ is isomorphic to $\tau$,
i.e.~$\sigma$ is a relabeling of $\tau$. There is an induced isomorphism
$\Psi: \mbar(X,\tau) \to \mbar(X,\sigma)$.

\begin{proposition}
\[
  \Psi^* \vso_{\mbar(X,\sigma)} = \vso_{\mbar(X,\tau)}, \quad
  \Psi_* \vso_{\mbar(X,\tau)} = \vso_{\mbar(X,\sigma)}.
\]
\end{proposition}

The proof is obvious.

\subsection{Contractions} \label{ss:3.7}
Let $\phi: \tau \to \sigma$ be a map of stable modular graphs by contracting
one edge or one loop, such that $genus(\sigma)=genus(\tau)$.
Denote $e$ the edge or loop in question (in $\tau$).
Let $\tilde{\sigma}$ be a fixed modular graph obtained from
$\sigma$ by adding $k$ tails and $(\tilde{\sigma} ,\beta)$ is a
stable modular graph with degree.

$(\tau_1^i, \beta_1^{i_j})$ are the stable modular graphs with degrees of the
following form.
$\tau_1^i$ are obtained from $\tau$ by adding $k$ tails in ways compatible
with $\tilde{\sigma} \to \sigma$. Thus there is a natural map from
$\tau_1^i$ to $\tilde{\sigma}$ by contracting $e$, and the following diagram
is commutative:
\[
 \begin{CD}
  \tau_1^i @>>> \tilde{\sigma} \\
  @VVV      @VVV \\
  \tau @>>>     \sigma
 \end{CD}
\]
where the two horizontal arrows contracting $e$ and vertical arrows are
forgetful maps (by forgetting $k$ new tails).
$\beta_1^{i_j}$ are the degrees on $\tau_1^i$, compatible with $\beta$.

Similarly, $(\tau_m^i, \beta_m^{i_j})$ are the stable modular graphs with
degrees such that $\tau_m^i$ is obtained from $\tau$ by first replacing
$e$ with a chain $c_m$ of $m$ edges and $m-1$ vertices such that the
genera of these $m-1$ vertices are all zero. One then adds $k$ tails to
this new graph in ways compatible with $\tilde{\sigma} \to \sigma$ and
contraction map $\tau_m^i \to \tilde{\sigma}$, which contracts $c_m$.
We again arrive at the following commutative diagram:
\[
 \begin{CD}
  \tau_m^i @>p_m>> \tilde{\sigma} \\
  @V{a_m}VV @VVV \\
  \tau @>>>     \sigma
 \end{CD}
\]
where $a_m$ is the composition of forgetting $k$ tails and contraction.
The degrees $\beta_m^{i_j}$ are the degrees assigned on $\tau_m^i$ in ways
compatible with $\beta$ under $p_m$.

Thus there is a natural morphism
$\phi_m: \coprod \mbar(X,\tau_m^i,\beta_m^{i_j}) \to
\mbar(X,\tilde{\sigma},\beta)$ and the diagram
\begin{equation} \label{e:contractions}
 \begin{CD}
  \coprod \mbar(X,\tau_m^i,\beta_m^{i_j}) @>{\phi_m}>>
    \mbar(X,\tilde{\sigma},\beta) \\
  @VV{\alpha_m}V    @VV{l}V \\
  \mbar_{\tau} @>{\Phi}>> \mbar_{\sigma}
 \end{CD}
\end{equation}
is commutative.

For example, if $\mbar_{\tau}=\mbar_{0,4} \times \mbar_{0,5}$ and
$\mbar_{\sigma}=\mbar_{0,7}$. Then
$\mbar(X, \tilde{\sigma}, \beta) = \mbar_{0,7+k}(X, \beta)$
and $\mbar(X, \tau_2, \beta_2)$ are of the form
$\mbar_{0,4+k_a}(X,\beta^a) \times_X \mbar_{0,2+k_b}(X, \beta^b) \times_X
\mbar_{0,5+k_c}(X,\beta^c)$ such that $k_a+k_b+k_c=k$ and
$\beta^a+\beta^b+\beta^c=\beta_2$.
Here the notation $\mbar_{g_1,n_1+1}(X,\beta_1) \times_X
\mbar_{g_2,n_2+1}(X,\beta_2)$ is the fibre product
\[
 \begin{CD}
  \mbar_{g_1,n_1+1}(X,\beta_1) \times_X \mbar_{g_2,n_2+1}(X,\beta_2) @>>>
  \mbar_{g_1,n_1+1}(X,\beta_1) \times \mbar_{g_2,n_2+1}(X,\beta_2) \\
  @VVV          @VV{\ev_{n_1+1} \times \ev_{1}}V \\
  X @>{\Delta}>> X \times X .
 \end{CD}
\]

The commutative diagram \eqref{e:contractions} induces a morphism
\begin{equation} \label{e:8}
  \Psi_{m}: \coprod \mbar(X,\tau_{m}^i,\beta_m^{i_j}) \to
  \mbar_{\tau} \times_{\mbar_{\sigma}} \mbar(X, \tilde{\sigma}, \beta).
\end{equation}

\begin{proposition} \label{p:contractions}
\begin{equation} \label{e:9}
 \sum_m (-1)^{m+1}  {\Psi_m}_*
    \sum_{(i, i_j)} \vso_{\mbar(X,\tau_m^i, \beta^{i_j}_m)}
  = \Phi^! (\vso_{\mbar(X,\tilde{\sigma},\beta)})
\end{equation}
in $K_{\circ}(\mbar_{\tau} \times_{\mbar_{\sigma}}
\mbar(X, \tilde{\sigma},\beta))$. The operation $\Phi^!$ on the RHS of
\eqref{e:9} is the refined Gysin map from
$K_{\circ}(\mbar(X,\tilde{\sigma},\beta))$ to $K_{\circ}(\mbar_{\tau}
\times_{\mbar_{\sigma}} \mbar(X,\tilde{\sigma},\beta))$.
\end{proposition}

\begin{proof}
Consider the following commutative diagram
\begin{equation*}
 \begin{CD}
  \coprod \mathfrak{M}_{\tau^i_m} @>>> \mathfrak{M}_{\tilde{\sigma}} \\
  @VVV      @VVjV \\
  \mbar_{\tau} @>\Phi>> \mbar_{\sigma}.
 \end{CD}
\end{equation*}
It induces a morphism
\[
 \psi_m : \coprod \mathfrak{M}_{\tau^i_m} \to
   \mbar_{\tau} \times_{\mbar_{\sigma}} \mathfrak{M}_{\tilde{\sigma}}.
\]
We will show an analogue of \eqref{e:9} in this setting. Namely,
\begin{equation} \label{e:10}
 \sum_{m,i} (-1)^{m+1}  {\psi_m}_* \so_{\mathfrak{M}_{\tau_m^i}}
  = \Phi^! \so_{\mathfrak{M}_{\tilde{\sigma}}}.
\end{equation}
Note that the vertical maps only contracts rational components with two marked
points (or less).
That is the reason why only chains of rational curves are relevant.

It is easy to see that
$\mbar_{\tau} \times_{\mbar_{\sigma}} \mathfrak{M}_{\tilde{\sigma}} \to
\mathfrak{M}_{\tilde{\sigma}}$ is a divisor $D_{\tau}$ of normal crossing in
$\mathfrak{M}_{\tilde{\sigma}}$, and the ``smoothing map''
\[
 \mu_1: \coprod \mathfrak{M}_{\tau_1^i} \to D_{\tau}
\]
is a finite, unramified and surjective birational morphism. In other words,
the smooth Artin stack $\coprod \mathfrak{M}_{\tau_1^i}$
separates the normal crossing point. Although $D_{\tau}$ is not equal to
$\mu_* \so_{\coprod \mathfrak{M}_{\tau_1^i}}$,
the difference is supported on the normal crossing subset. The following lemma
shows how one could obtain the structure sheaf
$\so_{D_{\tau}}$ by the inclusion-exclusion principle.

\begin{lemma}
Let $D=\cup_{i=1}^k D_i$ be a divisor with normal crossing, such that $D_i$
are smooth, disjoint away from the origin. Furthermore locally at origin
$D$ is defined by $x_1 \ldots x_k=0$. Then
\begin{equation} \label{e:iep}
  0 \to {\so}_D \to \sum_i \so_{D_i}  \to \sum_{i<j} \so_{D_i \cap D_j}
  \to \ldots \to \so_{D_1 \cap \ldots \cap D_k} \to 0
\end{equation}
is an exact sequence.
\end{lemma}

\begin{proof}
Equation \eqref{e:iep} is equivalent to the study the exactness of
the following sequence locally at the origin
\begin{multline*}
  0 \to {\so}/(x_1 x_2\ldots x_k) \to \sum_i {\so}/(x_i) \to \\
  \sum_{i<j} \so /(x_i, x_j) \to \ldots \to \so/(x_1, x_2, \ldots, x_k) \to 0.
\end{multline*}
This equation holds because, e.g.~$k=2$, the following sequence
\[
  0 \to {\so}/{(x_1 x_2)} \to {\so}/{(x_1)} \oplus {\so}/{(x_2)}
  \to {\so}/{(x_1, x_2)} \to 0
\]
is exact. In the general case this follows from the inclusion-exclusion
principle. For the exactness of \eqref{e:iep} at points away from the origin,
where $l$ divisors intersect ($l < k$), it follows from the induction as
the divisor locally defined by the equation $y_{i_1} \ldots y_{i_l}$.
\end{proof}

It remains to note that the normal crossing substack is stratified by the
image of
\[
  \mu_m : \coprod \mathfrak{M}_{\tau_m^i} \to D_{\tau}
\]
with codimension $m$ in $\mathfrak{M}_{\tilde{\sigma}}$. This concludes the
proof of \eqref{e:10}.

Now consider the following diagram
\[
 \begin{CD}
  \mbar_{\tau_m^i} (X, \beta_m^{i_j}) @>>> \mbar_{\tilde{\sigma}}(X, \beta) \\
  @VVV     @VVV \\
  \mathfrak{M}_{\tau_m^i} @>w>> \mathfrak{M}_{\tilde{\sigma}}
 \end{CD}
\]
which defines a compatible perfect obstruction theory over $w$. Therefore
by Proposition~\ref{p:functoriality}
\begin{equation} \label{e:11}
  w^! \vso_{\mbar_{\tilde{\sigma}}(X, \beta)}
 = \vso_{\mbar_{\tau_m^i} (X, \beta_m^{i_j})}.
\end{equation}
Proposition follows from the combination of equations \eqref{e:10} and
\eqref{e:11}.
\end{proof}

\section{Quantum $K$-invariants}

\subsection{Preliminaries on topological $K$-theory}
The $K$-theory used in this article has four variants:
\begin{itemize}
 \item $K_{\circ}(V)$: the Grothendieck group of coherent sheaves on $V$.
 \item $K^{\circ}(V)$: the Grothendieck group of algebraic vector bundles
    on $V$.
 \item $K^*(V):=K^0(V) \oplus K^1(V)$: the topological $K$-cohomology theory
    on $V$ (of complex vector bundles).
 \item $K_*(V):=K_0(V) \oplus K_1(V)$: the topological $K$-homology theory on
    $V$.
\end{itemize}

For the benefit of algebraic geometers let us recall some definitions and
useful properties of topological $K$-theory.
All the spaces involved are topological spaces or algebraic schemes. See
Remark~\ref{r:6} below for applications to stacks.

Let $V \hookrightarrow \cc^N$ be a topological closed embedding. Define
\[
  K_i(V) := K^{i}(\cc^N, \cc^N -V) \overset{\text{def}}{=} K^i_V(\cc^N).
\]
(Remember that there is a Bott periodicity.)
Note that this definition is independent of the closed embedding.

There is a \emph{cap product} between homology and cohomology theories. In
algebraic $K$-theory, it is just the tensor product. In topological
$K$-theory, it is defined as follows. Let $V \subset \cc^N$ be the closed
embedding as above and $U$ be a neighborhood of $V$ in $\cc^N$ such that
every element in $K^*(V)$ extends to an element in $K^* (U)$. $K_i(V)=
K^{i}(\cc^N, \cc^N-V) =K^{i}(U,U-V)$ by excision. Then the cap product is
defined as

\begin{multline*}
  K^i(V) \otimes K_j(V) \to K^i (U) \otimes K^{j}(U, U-V) \\
 \to K^{i+j}(U, U-V) = K_{j+i} (V) \overset{\cong}{\to} K_{j-i}(V),
\end{multline*}
where the second arrow is the cup product and the last arrow is the Bott
periodicity.

One can define the \emph{push-forward}
\[
  f_* : K_i(V') \to K_i(V)
\]
for $f:V' \to V$ a continuous proper map. Since $f$ is proper, there is
$\phi':V' \to D^m$ ($D^m$ is a polydisk in $\cc^m$) such that
\[ (f,\phi'): V' \to V \times D^m \]
is a closed embedding. Choose a closed
embedding $\phi: V \to \cc^N$ as above. The push-forward is the composition
\begin{multline*}
  K_i(V') \cong K^{i}(\cc^{N+m}, \cc^{N+m} - (\phi f, \phi')(X)) \\
  \overset{\text{res}}{\to} K^{i}(\cc^{N+m}, \cc^{N}\times
  (\cc^m - D^m) \cup (\cc^N -\phi(V)) \times \cc^m) \\
  \overset{\cong}{\leftarrow} K^{i} (\cc^N, \cc^N - \phi(V)) \cong K_i(V)
\end{multline*}
where the $\leftarrow$ is the Thom isomorphism, i.e.~multiplying the canonical
element of $K^{0}(\cc^m, \cc^m - D^m)$.

There is an obvious homomorphism $K^{\circ}(V) \to K^0(V)$ as
every algebraic vector bundle is a topological one.
One can also construct a homomorphism
\begin{equation} \label{e:14}
  Top: K_{\circ}(V) \to K_0(V)
\end{equation}
as follows. Choose a closed embedding of $V$ to a smooth scheme $i: V
\hookrightarrow Y$. (In our applications, all schemes are quasi-projective.)
For a bounded complex of coherent sheaves $\alpha^{\bullet}$ on
$V$, $i_*(\alpha^{\bullet})$ is quasi-isomorphic to a bounded complex
$E^{\bullet}$ of locally free sheaves on $Y$, which is exact out of $V$.
Regard $E^{\bullet}$ as a a topological element and embed $Y$ in some
$\cc^n$. By Thom isomorphism $E^{\bullet}$ is identified as a bounded complex
of topological vector bundles on $\cc^n$, exact off $V$, which is then
an element in $K_0(V)$.

A relation between the push-forward maps in two $K$-homology
theories is given by the following theorem of
Baum--Fulton--MacPherson:

\begin{theorem} \label{t:bfm} {\rm(\cite{BFM})}
For a proper morphism $g:X \to Y$, the following diagram:
\[
 \begin{CD}
  K_{\circ}(X) @>>> K_0(X) \\
  @V{g_*}VV @V{g_*}VV  \\
  K_{\circ}(Y) @>>> K_0(Y)
 \end{CD}
\]
commutes.
\end{theorem}

There is again a bivariant topological $K$-theory which unifies
the above setting. The homomorphism \eqref{e:14} is also valid in the
bivariant setting \cite{FM}.

\begin{remark} \label{r:6}
Although the $K$-theory of a moduli stack is quite different from
the $K$-theory of its coarse moduli space, one can define the quantum
$K$-invariants via coarse moduli space in the following way.
Let $p: \mathcal{M} \to M$ be the canonical map from the moduli stack to
its coarse moduli space. First construct the virtual structure sheaf
$\vso_{\mathcal{M}}$ on the moduli stack, then push-forward $p_*(\vso)$ to
its coarse moduli space, which is a topological space. One gets
\[
  \chi(\mathcal{M}, \ev_{\mathcal{M}}^*(\gamma)\, E \, \vso_{\mathcal{M}})
  = \chi(M, \ev_{M}^*(\gamma) \, p_*(E\, \vso_{\mathcal{M}}))
\]
by projection formula as $\ev:\mathcal{M} \to X$ factors through $M$.
$E$ is any $K$-element on $\mathcal{M}$, like $\cl$ or $\h$ (Hodge bundle).
This means that one gets the same quantum $K$-invariants on
coarse moduli space if one takes the suitable virtual structure sheaf.

Therefore, one can define quantum $K$-invariants for topological
$K$-theory on \emph{topological spaces} rather than orbispaces,
by using the topological $K$-homology for topological
space associated to coarse moduli space and $Top(p_*(\vso))$ as the virtual
structure sheaf. This eliminates the difficulty of constructing topological
$K$-homology on orbispaces, which will be addressed in another paper.
\end{remark}

\subsection{Definition of Quantum $K$-invariants} \label{ss:4.2}
In this subsection we will propose a construction of
$K$-theoretic invariants \footnote{The word `invariant' should be taken as
self-contained since we have not proved that it is a, say, symplectic
invariant as in the case of cohomology theory. This will be discussed
in a separate paper.} of the
Gromov--Witten type, which we call quantum $K$-invariants.

Let $K^{\bullet}(V)$ be the (algebraic or topological) $K$-cohomology of $V$
and $K_{\bullet}(V)$ be the (algebraic or topological) $K$-homology of $V$.
When $V$ is smooth, the notation $K(V)$ is used to denote both
$K^{\bullet}(V)$ and $K_{\bullet}(V)$ as they are isomorphic.

The quantum $K$-invariants is defined to be
\[
 \la \gamma_1, \gamma_2, \ldots, \gamma_n; F \ra_{g,n,\beta}:=
  \chi \left(\mbar_{g,n}(X,\beta), \vso \otimes
  \ev^*(\gamma_1 \otimes \ldots \gamma_n) \otimes \st^*(F) \right),
\]
where $\gamma_i \in K(X), F \in K(\mbar_{g,n})$ and $\chi$ is the push-forward
to a point ($\spec \cc$). Note that the triple $(g,n, \beta)$ is chosen so that
$\mbar_{g,n}(X,\beta)$ is defined. Note that due to Baum--Fulton--MacPherson
Theorem~\ref{t:bfm} topological invariants will be equal to algebraic
invariants (whenever applicable). We will not make distinctions.

One could also include the gravitational descendents and define:
\begin{equation}
 \begin{split}
   &\la \tau_{k_1}(\gamma_1), \tau_{k_2}(\gamma_2), \ldots,
   \tau_{k_n}(\gamma_n); \st^*(F) \ra_{g,n,\beta} \\
  :=&\chi \left( \mbar_{g,n}(X,\beta), \so^{\vir} \otimes
  \bigl(\otimes_{i=1}^n \cl_i^{\otimes k_i} \ev_i^*(\gamma_i) \bigr)
  \otimes \st^*(F) \right).
 \end{split}
\end{equation}

\emph{From this point we assume the $K$-theory is the topological
$K$-theory.} The modification to algebraic $K$-theory is in most places
straightforward. A key point will be discussed in Remark~\ref{r:9}.

Let $X$ be a smooth projective variety, and $E \subset H_2(X, \zz)$ denote 
the semigroup of effective curve classes. Let $\cc[E]$ be the semigroup ring
determined by $E$. Since $0 \in E$, $\cc[E]$ has a unit element. For
$\beta \in E$, the corresponding element of $\cc[E]$ will be denoted by
$Q^{\beta}$.

Let $\mathfrak{m}$ be the maximal ideal in $\cc[E]$ generated by the nonzero
elements of $E$. The \emph{Novikov ring} $N(X)$ is defined to be the 
completion of $\cc[E]$ in the $\mathfrak{m}$-adic topology. Alternatively,
$N(X)$ may be defined by formal series in $Q^{\beta}$:
\[
  N(X) = \{ \sum_{\beta \in E} c_{\beta} Q^{\beta} | c_{\beta} \in \cc \} .
\]

Let $e_0=\so, e_1, e_2, \ldots$ be a basis of $K(X)_{\qq}$ and let
$t_i$ be its dual coordinates.
Define $t:= \sum_i t_i e_i$. The $K$-theoretic Poincar\'e pairing
\[
   ( \cdot, \cdot ) : K(X)_{\qq} \otimes K(X)_{\qq} \to \qq
\]
is defined to be $( e_i, e_j) = \chi(e_i \otimes e_j)$. This is a perfect
pairing whenever $X$ is smooth. This defines a metric
$g_{ij}:=( e_i, e_j )$ on $K(X)$.

The \emph{quantum $K$-potential} of genus $0$ is a generating series of
genus zero quantum $K$-invariants:
\begin{equation} \label{e:gwp}
 G(t, Q) := \frac{1}{2}(t,t) + \sum_{n=0}^{\infty} \sum_{\beta \in E}
  \frac{Q^{\beta}}{n!} \la t,\ldots,t \ra_{0,n,\beta}.
\end{equation}
Here $(0,n,\beta)$ is a triple such that $\mbar_{0,n}(X,\beta)$ exists.
It is obvious that $G(t,Q) \in N(X) \otimes_{\cc} \cc [[t]]$.

We will see later that the metric in quantum $K$-theory ought to be
``quantized'' due to the modified contraction axiom (\S~3.7). Define the
\emph{quantum $K$-metric} to be
\begin{equation} \label{e:qm}
  ((e_i, e_j)) = G_{ij}:= \p_{t_i} \p_{t_j} G(t).
\end{equation}
$G_{ij}$ is the quantization of $g_{ij}$ because $G_{ij}|_{q=0} = g_{ij}$
by definition. Define $G^{ij}$ to be the inverse matrix of $G_{ij}$.

One could also define the quantum $K$-classes for $X$ to be a family of linear
maps (cf.~\cite{KM})
\[
  I^{X}_{g,n,\beta} : K(X)^{\otimes n} \to K(\mbar_{g,n})
\]
defined by
\begin{equation} \label{e:km}
  I^X_{g,n,\beta} (u_1 \otimes u_2 \otimes \ldots \otimes u_n)
 := \st_* (\so^{\vir}_{\mbar_{g,n}(X,d)} \otimes
           \ev_1^*(u_1) \otimes \ldots \otimes \ev_n^* (u_n)).
\end{equation}
The equivalence of these two definitions follows from the fact that
the Poincar\'e duality in $K(\mbar_{g,n},\qq)$ is a perfect pairing.

\begin{remark} \label{D:kgw}
1. It is obvious how to include the modular graphs into our data and define
quantum $K$-invariants associated to a modular graph.

2. Finally we mention that there is an \emph{equivariant} version of
these invariants. If there is an algebraic group (or compact Lie group) $G$
acting on $X$, then the moduli space $\mbar_{g,n}(X,\beta)$ (and its virtual
structure sheaf) has a $G$ action too.
We can talk about the $G$-equivariant quantum $K$-invariants, completely
parallel to equivariant quantum cohomology theory.
\end{remark}

\subsection{On Kontsevich--Manin axioms in $K$-theory}
The quantum $K$-invariants (without descendents) as defined above
obviously can not satisfy all axioms of cohomological Gromov--Witten
invariants \cite{KM}.
The \emph{effectivity} and \emph{motivic} axioms follow from the construction.
The \emph{grading} and \emph{divisor} axioms are missing, due to the lack of
``dimension counting''.
We will discuss briefly the remaining five axioms below.

\textbf{$S_n$-covariance:}
The quantum $K$-class is covariant under $S_n$ action on the marked points.

\textbf{Fundamental class:} Let $e_0$ be the ($K$-class of)
structure sheaf of $X$ and $\pi_n: \mbar_{g,n} \to \mbar_{g,n-1}$. Then
\begin{equation}
  I^X_{g,n,\beta}(u_1\otimes u_2 \otimes \ldots \otimes u_{n-1}\otimes
  e_0) = \pi_n^* (I^X_{g,n-1,\beta}(u_1\otimes u_2 \otimes \ldots
 \otimes u_{n-1})).
\end{equation}

\textbf{Mapping to a point:}
Suppose that $\beta=0$. Consider the following diagram
\[
 \begin{CD}
 \mbar_{g,n} \times X @>p_2>> X \\
  @VV{p_1}V \\
 \mbar_{g,n}
 \end{CD}
\]
where $p_i$ are the projection morphisms. Then
\begin{equation}
  I^X_{g,n,0} (u_1\otimes \ldots \otimes u_n)=
  (p_1)_* p_2^*(u_1\otimes \ldots \otimes u_n \otimes
   \lambda_{-1} (R^1\pi_* (\so_C \otimes T_V))^{\vee})
\end{equation}

\textbf{Splitting:}
Fix $g_1,g_2$ and $n_1,n_2$ such that $g=g_1+g_2, n=n_1+n_2$.
Let
\[
 \Phi: \mbar_{g_1,n_1+1} \times \mbar_{g_2,n_2+1} \to \mbar_{g,n}
\]
be the contraction map which glues the last marked point of
$\mbar_{g_1,n_1+1}$ to the first marked point of $\mbar_{g_2,n_2+1}$
(contracting a non-looping edge). Then
\begin{equation} \label{e:splitting}
 \begin{split}
  &\Phi^* \sum_{k,\beta} q^{\beta} \frac{1}{k!}
  \ft^k_* I^X_{g,n+k,\beta}(\gamma_1 \otimes \ldots \otimes \gamma_n \otimes
    t \otimes \ldots \otimes t) \\
  = &\sum_{ij}
  \left( \sum_{k_1, \beta'} q^{\beta'}
  \frac{1}{k_1 !} \ft^{k_1}_* I^X_{g_1,n_1+k_1+1,\beta'}
  (\gamma_1 \otimes \ldots \gamma_{n_1} \otimes
        t \otimes \ldots \otimes t \otimes e_i) \right) \\
  &G^{ij}(t)  \left(\sum_{k_2, \beta''}
  \frac{1}{k_2 !} \ft^{k_2}_* I^X_{g_2,n_2+k_2+1,\beta''}
  (e_j \otimes \gamma_{n_1+1} \otimes \ldots \gamma_n \otimes
        t \otimes \ldots \otimes t) \right),
 \end{split}
\end{equation}
where the notation $\ft^k: \mbar_{g,n+k} \to \mbar_{g,n}$ stands for the
forgetful map which forget the additional $k$ marked points.

\textbf{Genus reduction:}
Let
\[
 \Phi: \mbar_{g-1,n+2} \to \mbar_{g,n}
\]
be the contraction map which glues the last two marked points
(contracting a loop). Then
\[
 \begin{split}
  &\Phi^* \sum_{k,\beta} q^{\beta} \frac{1}{k!}
  \ft^k_* I^X_{g,n+k,\beta}(\gamma_1 \otimes \ldots \otimes \gamma_n \otimes
        t \otimes \ldots \otimes t) \\
 =&\sum_{ij} \left(\sum_{k, \beta} q^{\beta}  \frac{1}{k !}
  \ft^k_* I^X_{g+1,n+k+2,\beta}
 (\gamma_1 \otimes \ldots \gamma_{n} \otimes t \otimes \ldots \otimes
  t\otimes e_i \otimes e_j)\right)   G^{ij}(t)
 \end{split}
\]

The proofs of these five axioms follows from the corresponding axioms of
virtual structure sheaves. For example, Contractions, Cutting Edges and
Products imply the splitting and genus reduction axioms. The only novelty
is the appearance of $G^{ij}$ in Splitting and genus reduction axioms.
This has its origin in Contractions Axiom (Section~\ref{ss:3.7}).
Let us write $G_{ij}(t,q) = g_{ij} + F_{ij}(t,q)$.
The inverse matrix $G^{ij}(t,q)$ is therefore
\[
  g^{ij} + \sum_{m \ge 2} (-1)^{m+1} g^{i a_1} F_{a_1 b_1} g^{b_1 a_2} \ldots
  F_{a_{m-1} b_{m-1}} g^{b_{m-1} j}
\]
due to the matrix geometric series
\[
 \frac{1}{1+f} = 1 -f +f^2 -f^3 + \ldots.
\]
Note that the contributions from the contracted chains (denoted $c_m$)
in Proposition~\ref{p:contractions} are exactly
$F_{a_1 b_1} g^{b_1 a_2} \ldots F_{a_{m-1} b_{m-1}}$.
The fact that the combinatorics involved can be simplified by introducing
$G^{ij}$ is first observed in \cite{AG4}.

\subsection{String equation}  
Let $\pi: \mbar_{g,n+1}(X,\beta) \to
\mbar_{g,n}(X,\beta)$ be the forgetful map by forgetting the last marked point.
\footnote{It is proved in \cite{BM} that 
$\ft_{n+1} : \mbar_{g,n+1}(X,\beta) \to \mbar_{g,n}(X,\beta)$ is the universal
curve $\pi: \mathcal{C} \to \mbar_{g,n}(X,\beta)$. It is therefore 
convenient not to distinguish these two.}
String equation asks the relations between $\la \tau_{k_1}(\gamma_1), \ldots,
\tau_{k_n}(\gamma_n), \tau_0(1)\ra_{g,n+1,\beta}$ and
$\la \tau_{k_1}(\gamma_1), \ldots, \tau_{k_n}(\gamma_n) \ra_{g,n,\beta}$.
In the higher genus case, the Hodge bundle is also included.
It is easy to see that the relation is reduced to the following equations
(by projection formula):

For $g=0$:
\begin{equation} \label{e:string0}
 \pi_*\Bigl(\so^{\vir}\bigl(\prod_{i=1}^{n} \frac{1}{1-q_i \cl_i}\bigr)\Bigr)
 = \Bigl( 1+ \sum_{i=1}^{n} \frac{q_i}{1-q_i} \Bigr) \biggl( \so^{\vir}
   \Bigl(\prod_{i=1}^{n} \frac{1}{1-q_i \cl_i} \Bigr) \biggr) ,
\end{equation}
where $q_i$ are formal variables. The LHS and RHS are equal as formal series
in $q_i$'s.

For $g \ge 1$:
\begin{equation} \label{e:stringg}
 \begin{split}
   &\pi_* \left(\so^{\vir} \frac{1}{1-q\h^{*}}
    \prod_{i=1}^{n} \frac{1}{1-q_i \cl_i} \right) \\
 =  &\so^{\vir} \frac{1}{1-q \h^{*}} \Bigg[
   \left( 1-\h^{*}+\sum_{i=1}^{n} \frac{q_i}{1-q_i} \right)
  \left( \prod_{i=1}^{n}
  \frac{1}{1-q_i \cl_i} \right) \Bigg],
 \end{split}
\end{equation}
where $\h$ is the Hodge bundle,
i.e.~$\h:=R^0 \pi_* \omega_{\mathcal{C}/\mbar}$.

\begin{proof}
First notice that the virtual structure sheaves in the above equations can
be simultaneously erased from both sides by the Fundamental class and
Forgetting tails axioms.

First the case $g=0$. 
Let $D_i$ be the divisors on $\mbar_{0,n+1}(X,\beta)$ such that
the generic curves have two components: one contains the $i$-th and
$n+1$-th marked points, and the other contains the rest. It is well-known that
\begin{equation} \label{e:25}
  \so(\cl_i) = \pi^* l_i \otimes \so(D_i)
\end{equation}
(see e.g.~\cite{EW}), where for the notational convenience we have used
$l_i$ and $\cl_i$ for $i$-th universal cotangent line bundles on
$\mbar_{0,n}(X,\beta)$ and $\mbar_{0,n+1}(X, \beta)$ respectively.
It is easy to see that (\cite{YL1} Lemma~2)
\begin{equation} \label{e:lemma}
  R^1 \pi_* (\so(\sum_i d_i D_i))=0
\end{equation}
for $d_i \ge 0$. By projection formula
\[
  R^0 \pi_* (\otimes_{i=1}^n \cl^{\otimes d_i}_i)
  = ( \otimes_{i=1}^n l^{d_i}_i ) R^0 \pi_* \so(\sum_i d_i D_i).
\]
It remains to compute $H^0:= R^0 \pi_* \so(\sum_i d_i D_i)$, which is a
vector bundle because $R^1=0$.

To compute $H^0$, we use some local arguments.
The fibre $C_x$ of $\pi$ is a rational curve with $n$ marked
points $x_1, \ldots, x_n$. An element of $H^0(C_x,\so(\sum_i d_i D_i))$ is
a rational function with poles of order no more than $d_i$ at $x_i$.
Therefore $H^0(C_x, \so(\sum_i d_i D_i))$ is filtered by the degrees
of poles at $x_i$:
\[
  \mathcal{F}_0 \subset \mathcal{F}_1 \subset \ldots \subset
  \mathcal{F}_{d_i}.
\]
It is by definition that the graded piece $\mathcal{F}_{k+1}/
\mathcal{F}_k$ are isomorphic to $T^{\otimes k}_{x_i}$. Summarizing,
\begin{equation} \label{e:21}
  H^0 \overset{\text{$K$-theory}}{=} \so \oplus l^{-1}_1 \oplus \ldots
  l^{-d_1} \oplus l^{-1}_2 \oplus \ldots \oplus l^{-1}_n \oplus \ldots
  l^{-d_n}_n,
\end{equation}
where ``$K$-theory'' above means the equation is valid only in
$K$-theory as we have used the graded pieces.

It is now a matter of elementary computation (see Proposition~1
in \cite{YL1}) to deduce the genus zero string equation from \eqref{e:21}.

In the case $g \ge 1$ the Hodge bundles $\h$ will naturally occur.
Notice that $\h$ is actually equal to ${\pi}^*(\h)$, and this implies
that the factor $1/(1-q\h)$ commutes with $\pi_*$ due to the projection
formula. By Grothendieck--Riemann--Roch formula, \eqref{e:stringg} is a
rational function of $q$ and $q_i$'s. Therefore
\begin{equation} \label{e:a1}
\pi_* (\frac{1}{1-q\cl})= \pi_* (\frac{-q^{-1}\cl^{-1}}{1-q^{-1}\cl^{-1}}).
\end{equation}
Now for $d_i \geq 1$, we want to show
\begin{align*}
 &\pi_* \lp\otimes_{i=1}^{n} \cl_i^{-d_i}\rp \\
 = &\otimes_{i=1}^{n} \cl_i^{-d_i}
 \lp\tlb - \h^{*} -\sum_{i,\ d_i \neq 0}
 \lp\tlb + \cl_i+ \ldots +\cl_i^{d_i -1}\rp\rp .
\end{align*}
The above equality can be proved using the same arguments as in the
genus zero case. By \eqref{e:lemma} and projection formula
\[
  \pi_*(\otimes_{i=1}^{n} \cl_i^{-d_i})= (\otimes_{i=1}^{n} (l_i)^{-d_i})
   \, \pi_* \so(\sum_i -d_i D_i).
\]
Because $d_i \ge 1$ for all $i$, $R^0 \pi_*(\so(\sum_i -d_i D_i))=0$.
\[
  H^1:= R^1 \pi_* (\so(\sum_i -d_i D_i))
\]
is a vector bundle.

Now by Serre duality $(H^1)^* = R^0 \pi_* \omega(\sum_i d_i D_i)$, where
$\omega$ is the dualizing sheaf. Fibrewisely, $H^0(C, \omega(\sum_i d_i x_i))$
is the holomorphic differential with poles of order at most $d_i$ at $x_i$.
Thus we have a filtration
$\mathcal{F}_1 \subset \mathcal{F}_2 \subset \ldots \subset \mathcal{F}_{d_i}$
of degrees of poles at each marked point $x_i$ as above and the graded
bundles $\mathcal{F}_{k+1}/\mathcal{F}_k$ is isomorphic to
$l_i^{\otimes -k}$ and $\mathcal{F}_0=\h$. Therefore
\[
  \begin{split}
  &\pi_*(\otimes_{i=1}^{n} \cl_i^{-d_i}) \\
  =&\otimes_{i=1}^{n} (l_i)^{-d_i}
	\lp - R^1 \pi_* \so(\sum_i -d_i D_i) \rp  \\
  =& \otimes_{i=1}^{n} (l_i)^{-d_i}  
	\lp - (R^0 \pi_* \omega (\sum_i -d_i D_i))^* \rp  \\
  = &\otimes_{i=1}^{n} l_i^{-d_i} \lp\tlb - \h^{*} -\sum_{i,\ d_i \neq 0}
 	\lp \tlb + l_i+ \ldots +l_i^{d_i -1}\rp \rp.
  \end{split}
\]
Notice again that the last equality holds only in $K$-theory (using graded
objects). This means
\begin{align*}
 &\pi_* \lp \prod_{i=1}^{n} \frac{q_i^{-1}\cl_i^{-1}}{1-q_i^{-1}
    \cl_i^{-1}}\rp\\
  =&\lp 1-\h^{*}+\sum_{i=1}^{n} \frac{q_i^{-1}}{1-q_i^{-1}} \rp
 \lp \prod_{i=1}^{n} \frac{q_i^{-1} l_i^{-1}}{1-q_i^{-1} l_i^{-1}}\rp ,
\end{align*}
which is equivalent to \eqref{e:stringg}.
\end{proof}

\subsection{Dilaton equation}
Let $\pi: \mbar_{g,n+1}(X,\beta) \to \mbar_{g,n}(X,\beta)$ as above. The
dilaton equation asks the relations between $\la \tau_{k_1}(\gamma_1),\ldots,
\tau_{k_n}(\gamma_n),\tau_1(e_0) \ra_{g,n+1,\beta}$ and
$\la \tau_{k_1}(\gamma_1),\ldots, \tau_{k_n}(\gamma_n) \ra_{g,n,\beta}$.
\begin{equation} \label{e:dilaton}
 \begin{split}
   &\pi_* \left(\so^{\vir}  \frac{1}{1-q\h}  \biggl(
   \prod_{i=1}^{n} \frac{1}{1-q_i \cl_i} \biggr)  \cl_{n+1} \right) \\
 =  &\so^{\vir}  \frac{1}{1-q \h }  
   \Bigg[
   \left( \h -1 +\sum_{i=1}^{n} \frac{1}{1-q_i} 
  \right)   \left( \prod_{i=1}^{n} \frac{1}{1-q_i \cl_i} \right) \Bigg]
 \end{split}
\end{equation}

\begin{proof}
The proof is similar to the above proof of string equation.
In the proof we again use $l_i$ for $\cl_i$ on $\mbar_{g,n}(X,\beta)$.
We will use Forgetting tail axiom, \eqref{e:25}\eqref{e:lemma} and projection formula
\begin{align*}
  \pi_* (\otimes_{i=1}^n \cl_i^{d_i} \otimes \cl_{n+1})
 =&\otimes_{i=1}^n l_i^{d_i} \otimes
    \pi_* \left(\cl_{n+1} (\sum_i d_i D_i)\right) \\
 \intertext{Because of the relation $\cl_{n+1}=\omega(x_1+\ldots+x_n)$,}
 =&\otimes_{i=1}^n l_i^{d_i} \otimes  \pi_*
  \omega \Bigl( \sum_i (d_i+1) D_i \Bigr)\\
 =&\otimes_{i=1}^n l_i^{d_i} \otimes \Bigl( \h -1 +n + \sum_{i=1}^n
  \sum_{k=1}^{d_i} l_i^{-k} \Bigr).
\end{align*}
It is easy to see that this is equivalent to \eqref{e:dilaton}.
\end{proof}

\begin{remark}
The above computation also yields
\[
  \pi_* (\cl_{n+1} -1) = \h + \h^{*} +(n-2)1.
\]
This resembles the cohomological dilaton equation in the sense that the
rank of the RHS is $2g-2+n$.
\end{remark}

\begin{remark}
The techniques used in proving String and Dilaton equations can also be used
to find the push-forward of negative powers of $\cl$.
\end{remark}

\section{Quantum $K$-ring and Frobenius manifold}

In this section we generalize Givental's treatment \cite{AG4} to
non-convex algebraic manifolds. Many statements already appeared in
\cite{AG4} (albeit with the convexity condition there) and are included
here for completeness.

\subsection{The ring structure in quantum $K$-theory} \label{ss:5.1}
The definition of the ring structure of quantum $K$-theory is analogous to
that of quantum cohomology theory.

Let $\{e_i\}$ be a basis of $K(X)$ and $t_i$ its coordinates.
$t:= \sum_i t_i e_i$ as in Section~\ref{ss:4.2}.

%

The quantum $K$-product $*$ is defined to be:
\begin{equation} \label{E:41-2}
 ((e_i * e_j, e_k)) := \p_{t_i} \p_{t_j} \p_{t_k} G(t),
\end{equation}
where $G(t)$ is the genus zero potential defined in \eqref{e:gwp} and
$((\cdot,\cdot))$ is the quantized metric defined in equation \eqref{e:qm}.


The main result of this section is

\begin{theorem} \label{T:41-1}
$(K(X, \qq[[Q]]), *)$ is a commutative and
associative algebra, which deforms the usual ring structure of $K(X)$.
\end{theorem}

\begin{proof}
The deformation properties follow immediately from the
definition. Setting $q=0$ would make $*$ the ordinary tensor product.
The (super-)commutativity also follows from the definition.

The associativity follows from the Splitting Axiom \eqref{e:splitting}.
It is easy to see that the $K$-theoretic WDVV equation
\[
  \sum_{\mu \nu} G_{i j \mu} G^{\mu \nu} G_{\nu k l}
  = \sum_{\mu \nu} G_{i k \mu} G^{\mu \nu} G_{\nu j l}
  = \sum_{\mu \nu} G_{i l \mu} G^{\mu \nu} G_{\nu j k}
\]
is equivalent to the associativity of the quantum $K$-ring. The
WDVV equation follows from Splitting Axiom: Let $n=4, n_1=2, n_2=2$ in
\eqref{e:splitting}. Different ways of splitting amount to the same invariants.
\end{proof}

In the case $t=0$, we have defined a deformation by quantum three-point
function, which is similar to the \emph{pair-of-pants} structure. However,
in this case the metric ought to be quantized as well (and is replaced by
two-point function at $t=0$).

\begin{remark}
For convex $X$, a proof can be found in \cite{AG4}.
\end{remark}

\begin{remark} \label{r:9}
The above description is good for topological $K$-theory. For algebraic
$K$-theory, a finite basis does not necessarily exist and the
$K$-theoretic Poincar\'e pairing has a huge kernel. However, one could
proceed by using a presentation which is basis-free and does not rely on
Poincar\'e pairing. Namely, one could define the quantum $K$-product
(at $t=0$) as
\[
 e_i * e_j = \sum_{\beta} Q^{\beta} (\ev_3)_* \left( \ev_1^*(e_i) \ev_2^*(e_j)
 \sum_m (-1)^{m+1} \sum_{(i,i_j)} {\Psi_m}_*
    \vso_{\mbar(X,\tau^i_{m},\beta^{i_j}_m)} \right)
\]
where $\ev_i : \mbar_{0,3}(X,\beta) \to X$ are the evaluation morphisms and 
\[
  \sum_m (-1)^{m+1} \sum_{(i,i_j)} 
  {\Psi_m}_* \vso_{\mbar(X,\tau^i_{m},\beta^{i_j}_m)}
\]
is the alternating sum of virtual structure sheaves appeared in \eqref{e:9}.
\end{remark}

\subsection{Quantum $K$-theory and Frobenius (super-)manifolds}
In this section, we show that the construction of quantum $K$-ring produces
a (new) class of Frobenius manifolds, in the sense of \cite{YM}.
This in particular gives a positive answer to the question raised by
Bayer and Manin (\cite{BayerManin} 1.1.1 Question).
Note that in this case the identity element $e_0$ is \emph{not} flat.

The data of $K$-theoretic Frobenius structure includes:

(1) A metric $((e_i, e_j)) = G_{ij}:= \p_i \p_j G(t)$ on tangent spaces to
$K(X)_{\cc}$.

(2) The quantum multiplication with structural constants
\[ ((e_i * e_j,e_k))= \p_i \p_j \p_k G(t).
\]
on the tangent bundle.

(3) A connection on the tangent bundle defined by the operators of
quantum multiplication:
\[  \nabla_q:= d - \frac{1}{1-q} \sum_i (e_i *) d t_i.
\]

\begin{proposition}
1. The metric and the quantum multiplication $*$ define on the tangent
bundle a formal commutative associative Frobenius algebra with unit $1$.

2. The connection $\nabla_q$ are flat for $q \ne 1$.

3. The operator $\nabla_{-1}$ is the Levi-Civita connection of the metric
$((\cdot,\cdot))$.

4. The metric $((\cdot,\cdot))$ is flat.
\end{proposition}

\begin{proof} (see \cite{AG4})
The first two are formal consequences of WDVV equation. The Levi-Civita
connection of the metric $G_{ij}(t)$ is
\[
  \Gamma_{ij}^k = \frac{1}{2}(G_{il;j}+G_{jl;i}-G_{ij;l}) G^{lk}
  =\frac{1}{2} G_{ijl} G^{lk}= \frac{1}{2} (e_j)_i^k,
\]
where $G^{ij}$ is the inverse matrix of $G_{ij}$. Thus \emph{3.}~holds.
\emph{4.}~is an obvious consequence of \emph{3}.
\end{proof}

\begin{corollary}
$QK^*(X)_{\cc}$ is a formal Frobenius manifold over the Novikov ring $N(X)$.
\end{corollary}

\subsection{Quantum differential equation in $K$-theory}
Parallel to the discussion in quantum cohomology, we may introduce
quantum differential equation \cite{AG2}:
\[  \nabla_q S=0, \]
which is a system of linear partial differential equations.

\begin{theorem} \label{t:qde}
The matrix $(S_{ij})$
\[
  S_{ij}(t,Q):=g_{ij}+ \sum_{n,\beta} \frac{Q^{\beta}}{n!}
  {\Big\la} e_i, t,\ldots,t, \frac{e_j}{1-q \cl} {\Big\ra}_{0,n+2,\beta}
\]
is the fundamental solution to the $K$-theoretic quantum differential
equation. Namely, the column vectors form a complete set of solutions.
\end{theorem}

\begin{proof}
The proof relies on two ingredients: the string equation and WDVV equation.
However, the WDVV equation here takes a slightly generalized form. Put
$e_i,e_j,e_k, \frac{e_l}{1-q \cl}$ on the distinguished four marked points.
The proof of WDVV equation goes through and we get
\[
  G_{i,j,\alpha} G^{\alpha \beta} \p_k S_{\beta l}
  = G_{i,k,\alpha} G^{\alpha \beta} \p_j S_{\beta l},
\]
or equivalently,
\begin{equation} \label{e:gwdvv}
  (e_j *) \p_k S = (e_k *) \p_j S.
\end{equation}
Put $j=0$. Since $(e_0 *)=\operatorname{Id}$, \eqref{e:gwdvv}
implies that
\[
  \p_k S = (e_k *) \p_0 S.
\]
Now by string equation
\[
   \p_0 S = \frac{1}{1-q} S.
\]
\end{proof}


\begin{thebibliography}{40}

\bibitem{sga6}
SGA 6, Dirig\'e par P. Berthelot, A. Grothendieck et L. Illusie.
Lecture Notes in Mathematics, Vol. \textbf{225}, Springer-Verlag,
Berlin-New York, 1971.

\bibitem{BFM}
P.~Baum, W.~Fulton, R.~MacPherson, \emph{Riemann-Roch and topological $K$
theory for singular varieties}, Acta Math. \textbf{143} (1979), no. 3-4,
155--192.

\bibitem{BayerManin}
A.~Bayer, Yu.~Manin, \emph{(Semi)simple exercises in quantum cohomology},
math.AG/0103164.

\bibitem{KB}
K.~Behrend, \emph{Gromov--Witten invariants in algebraic geometry}, Invent.
Math. \textbf{127} (1997), no. 3, 601-617.

\bibitem{BF}
K.~Behrend, B.~Fantechi, \emph{The intrinsic normal cone}, Invent. math.
\textbf{128} (1997), no. 1, 45-88.

\bibitem{BM}
K.~Behrend, Yu.~Manin, \emph{Stacks of stable maps and Gromov--Witten
invariants}, Duke math. J. \textbf{85} (1996), no. 1, 1-60.


\bibitem{WF}
W.~Fulton, \emph{Intersection theory}, Springer-Verlag, Berlin, 1998.

\bibitem{FL}
W.~Fulton, S.~Lang, \emph{Riemann--Roch algebra},
Grundlehren der Mathematischen Wissenschaften
[Fundamental Principles of Mathematical Science], 277.
Springer-Verlag, New York-Berlin, 1985. x+203 pp

\bibitem{FM}
W.~Fulton, R.~MacPherson, \emph{Categorical framework for the study of
singular spaces}, Mem. Amer. Math. Soc. \textbf{31} (1981), no. 243, vi+165 pp.

\bibitem{FP}
W.~Fulton, R.~Pandharipande, \emph{Notes on stable maps and quantum
cohomology}, Algebraic geometry---Santa Cruz 1995, Proce. Symp. Pure. Math.
62. Part 2, (1997) 45-96.

\bibitem{AG2}
A.~Givental, \emph{Equivariant Gromov-Witten invariants},
IMRN \textbf{13} (1996) 613-63.

\bibitem{AG4}
A.~Givental, \emph{On the WDVV-equation in quantum $K$-theory},
Mich. Math. J. \textbf{48} (2000) 295-304.

\bibitem{GK}
A.~Givental, B.~Kim, \emph{Quantum cohomology of flag manifolds and Toda
lattices}, Comm. Math. Phys. 168 (1995), no. 3, 609--641.

\bibitem{GL}
A.~Givental, Y.-P.~Lee, \emph{Quantum $K$-theory on flag manifolds,
finite-difference Toda lattices and quantum groups}, math.AG/0108105, to appear
in Invent. Math.

\bibitem{GP}
T.~Graber, R.~Pandharipande, \emph{Localization of virtual fundamental
class}, Invent. math. \textbf{135} (1999) 2 487-518.

\bibitem{BK}
B.~Kim, \emph{Quantum cohomology of flag manifolds $G/B$ and quantum Toda
lattices}, Ann. of Math. (2) 149 (1999), no. 1, 129--148.

\bibitem{MK1}
M.~Kontsevich, \emph{Intersection theory on the moduli space of curves and
the matrix Airy function}, Comm. Math. Phys. \textbf{147} (1992), no. 1, 1--23.

\bibitem{MK}
M.~Kontsevich,
\emph{Enumeration of rational curves via torus actions}, In: The
moduli space of curves. R.~Dijkgraff, C.~Faber, G.~van der Geer (Eds.),
Progress in Math., \textbf{129}, Birkh\"auser, 1995, 335-68.

\bibitem{KM}
M.~Kontsevich, Yu.~Manin, \emph{Gromov--Witten Classes, quantum cohomology,
and enumerative geometry}, Comm. Math. Phys. \textbf{164} (1994), no. 3,
525-562.

\bibitem{AK}
A.~Kresch, \emph{Canonical rational equivalence of intersections of divisors},
Invent. Math. \textbf{136} (1999), no. 3, 483--496.

\bibitem{YL1}
Y.-P.~Lee, \emph{A Formula for Euler characteristics of tautological line
bundles on the Deligne-Mumford moduli spaces}, IMRN, 1997, No. 8, 393-400.

\bibitem{YL2}
Y.-P.~Lee,
\emph{Euler characteristics of universal cotangent line bundles on
$\mbar_{0,1}$}, math.AG/0005217.

\bibitem{YL4}
Y.-P.~Lee, \emph{Quantum $K$-theory II: computations and open problems},
in preparation.

\bibitem{LP}
Y.-P.~Lee, R.~Pandharipande, \emph{A reconstruction theorem in quantum
cohomology and quantum $K$-theory}, preprint math.AG/0104084.

\bibitem{LT}
J.~Li, G.~Tian, \emph{Virtual moduli cycles and Gromov--Witten invariants},
J. Amer. Math. Soc. \textbf{11} (1998), no. 1, 119-174.

\bibitem{YM}
Yu.~Manin, \emph{Frobenius manifolds, quantum cohomology, and moduli spaces},
American Mathematical Society Colloquium Publications, \textbf{47},
American Mathematical Society, Providence, RI, 1999. 303 pp.

\bibitem{JM}
J.~Morava, \emph{Quantum generalized cohomology}, in \emph{Operads:
Proceedings of Renaissance Conferences} (Hartford, CT/Luminy, 1995), 407--419,
Contemp. Math., 202, Amer. Math. Soc., Providence, RI, 1997.

\bibitem{EW}  E.~Witten,
\emph{Two-dimensional gravity and intersection theory on moduli space},
Surveys in Diff. Geo. vol.~\textbf{1} (1991) 243-310.

\end{thebibliography}
\end{document}